\documentclass[11pt]{article}
\usepackage[T1]{fontenc}
\usepackage{lmodern}
\usepackage{amsmath,amssymb,amsthm,mathtools}
\usepackage[a4paper,margin=25mm]{geometry}
\usepackage{microtype}
\usepackage{xcolor}
\usepackage[colorlinks=true,linkcolor=blue,citecolor=blue,urlcolor=blue,
  hyperfootnotes=false,bookmarksnumbered]{hyperref}
\hypersetup{pdftitle={Ground sets of intersecting set-pair systems and maximal intersecting families},pdfauthor={Andrey Kupavskii}}
\allowdisplaybreaks[1]

\newtheorem{theorem}{Theorem}
\newtheorem{lemma}[theorem]{Lemma}
\newtheorem{corollary}[theorem]{Corollary}
\newtheorem{claim}{Claim}[section]
\theoremstyle{definition}
\newtheorem{construction}{Construction}

\newcommand{\A}{\mathcal A}
\newcommand{\B}{\mathcal B}
\newcommand{\F}{\mathcal F}

\newcommand{\PP}{\mathbb P}
\newcommand{\EE}{\mathbb E}

\title{Ground sets of intersecting set-pair systems\\
and maximal intersecting families}
\author{
Elizaveta Iarovikova\thanks{Moscow Institute of Physics and Technology;
E-mail: \url{liza.fm@yandex.ru}} \and
Andrey Kupavskii\thanks{Moscow Institute of Physics and Technology,
Saint-Petersburg State University, Innopolis University;
E-mail: \url{kupavskii@ya.ru}} \and
Hong Liu\thanks{Extremal Combinatorics and Probability Group,
Institute for Basic Science, Daejeon, South Korea;
E-mail: \url{hongliu@ibs.re.kr}} \and
Ilya Lobatskii\thanks{Moscow Institute of Physics and Technology;
E-mail: \url{lobzik22848@gmail.com}} \and
Yakov Shubin\thanks{Moscow Institute of Physics and Technology;
E-mail: \url{shubin.yakoff@gmail.com}}
}

\begin{document}
\maketitle

\begin{abstract}
An intersecting set-pair system is a family $\{(A_i,B_i):1\le i\le m\}$ such that
$A_i\cap B_j=\varnothing$ if and only if $i=j$. Tuza asked how large its ground set can be when $|A_i|=a$ and $|B_i|=b$ for every $i$. We determine this maximum asymptotically, proving that
\[
 n(a,b)=(1+o(1))\max_{1\le s\le a}s\binom{a+b-s}{b},
\]
for integers $a\ge b\ge1$ as $a\to\infty$. We also resolve asymptotically a question of Erd\H{o}s and Lov\'asz on the maximum size of the ground set of a maximal intersecting
$k$-uniform family. We prove that this maximum is $(1/2+o(1))\binom{2k-2}{k-1}$ as $k\to\infty$, showing that their construction is asymptotically optimal.
\end{abstract}

\section{Introduction}

We write $[m]=\{1,\ldots,m\}$. Given a set $X$ and an integer $k \geqslant 0$, denote by $2^X$ and $\binom{X}{k}$ the collections of all subsets and all $k$-element subsets ($k$-sets) of $X$, respectively. 

One classical object in extremal set theory is an intersecting set-pair system. A collection of pairs $\{(A_i,B_i):i\in[m]\}$ is called an
\emph{intersecting set-pair system} (ISP for short) if
\[
 A_i\cap B_j=\varnothing \text{ if and only if } i=j
 \text{ whenever } 1\le i,j\le m.
\]
Its \emph{ground set} is $V=\bigcup_{i=1}^m(A_i\cup B_i)$. If $|A_i|=a$ and $|B_i|=b$ for every $i\in[m]$, we call the system an $\operatorname{ISP}(a,b)$.

Intersecting set-pair systems arise naturally in a variety of extremal problems, including graph saturation~\cite{bollobas}, critical hypergraphs~\cite{tuza-critical}, and the enumeration of maximal intersecting families~\cite{nagy-patkos}.

An important property of ISPs is that we can bound their size independently of the ground set. Recall the following celebrated theorem of Bollob\'as.

\begin{theorem}[Bollob\'as~\cite{bollobas}]\label{bollobas}
Every ISP $\{(A_i,B_i):1\le i\le m\}$ satisfies
\begin{equation}\label{eq:bollobas}
 \sum_{i=1}^m\binom{|A_i|+|B_i|}{|A_i|}^{-1}\le1.
\end{equation}
\end{theorem}

This fundamental inequality has several proofs, including
combinatorial proofs~\cite{bollobas,jaeger-payan,katona-set-pairs}
and a linear-algebraic proof of the uniform version~\cite{lovasz-flats}.

Bollob\'as's inequality has numerous extensions, including skew versions~\cite{frankl-skew}, weakly cross-intersecting
systems~\cite{tuza-weak,kiraly-et-al}, and analogues for vector spaces~\cite{lovasz-flats}.
For more recent developments, see Scott and Wilmer~\cite{scott-wilmer} and Heged\"us and Frankl~\cite{hegedus-frankl}.

Theorem~\ref{bollobas} gives the sharp bound $m\le\binom{a+b}{b}$ on the number of pairs in an $\operatorname{ISP}(a,b)$. The following construction attains equality.

\begin{construction}\label{constr1}
Fix integers $a,b\ge1$, and let $X$ be a set of size $a+b$.
The pairs
\[
 \bigl\{(S,X\setminus S):S\in\tbinom{X}{a}\bigr\}
\]
form an $\operatorname{ISP}(a,b)$ with $m=\binom{a+b}{b}$ and $|V|=a+b$.
\end{construction}

Beyond the tightness of the bound in  Theorem~\ref{bollobas}, however, we have limited understanding of large ISPs. Motivated by that, Tuza~\cite{tuza-critical} suggested to look at the following parameter. Let $n(a,b)$ be the maximum size of the ground set of an $\operatorname{ISP}(a,b)$.  Construction~\ref{constr1} has a very small ground set. Tuza, however, noted that its modification  gives an ISP with a ground set that is exponential in the uniformity.

\begin{construction}\label{construct_ground}
Fix integers $a, b \ge 1$ and $1\le s\le a$. Let $X$ be a set of size $a+b-s$. For every $S\in\binom{X}{a-s}$, choose a set $P_S$ of size $s$, with all the sets $P_S$ pairwise disjoint and disjoint from $X$. Then
\[
 \bigl\{(S\cup P_S,X\setminus S):S\in\tbinom{X}{a-s}\bigr\}
\]
is an $\operatorname{ISP}(a,b)$ with $m=\binom{a+b-s}{b}$ and $|V|=s\binom{a+b-s}{b}+a+b-s$.
\end{construction}
Note that setting $s=0$ in Construction~\ref{construct_ground} gives Construction~\ref{constr1}. Put
\begin{equation}\label{eq:K}
 K(a,b)=\max_{1\le s\le a}s\binom{a+b-s}{b}.
\end{equation}
The maximum in~\eqref{eq:K} is attained at $s=\lceil a/(b+1)\rceil$.

Tuza conjectured that Construction~\ref{construct_ground} is best possible. He managed to prove that, for $a\ge b\ge1$,
\[
 \frac14\binom{a+b+1}{b+1}
 < n(a,b)
 \le
 1+\sum_{i=1}^{2b-2}\binom{i}{\lfloor i/2\rfloor}
  +\sum_{i=2b-1}^{a+b-1}\binom{i}{b}.
\]
In the case $a=b$, it gives $\left(\frac12+o(1)\right)\binom{2a}{a} \le n(a,a)\le(1+o(1))\binom{2a}{a}$ as $a\to\infty$. Our first result confirms the conjecture of Tuza asymptotically.

\begin{theorem}\label{thm:fixed}
For integers $a\ge b\ge1$ and $a\to\infty$, we have
\[
 n(a,b)=(1+o(1))K(a,b).
\]
\end{theorem}

In particular,  $n(a,a)\sim\frac12\binom{2a}{a}$. We should note that Nagy and Patk\'os \cite{nagy-patkos} corrected the conjecture of Tuza, by giving a slightly bigger construction for $a=b=k$ (the construction stays essentially the same, but  $(k-1,k-1)$-pairs replace the $(k,k-2)$-pairs inside the core $X$). It increases the lower bound on the size of the ground set by about $\frac 1k{2k-2\choose k-1}.$ 

Theorem~\ref{thm:fixed} also gives an asymptotically sharp bound on the sum of the number of pairs and the size of ground set. This yields a weak stability statement for Bollob\'as's inequality.

\begin{corollary}\label{cor:pairs-ground}
Every $\operatorname{ISP}(a,b)$ with $m$ pairs and ground set $V$
satisfies
\[
 m+|V|\le(1+o(1))
 \max_{0\le s\le a}(s+1)\binom{a+b-s}{b},
\]
for $a\ge b\ge1$ and $a\to\infty$. For $a=b$, this gives $m+|V|\le(1+o(1))\binom{2a}{a}.$ 
\end{corollary}

\begin{proof}
Adding a distinct new element to each $A_i$ gives an
$\operatorname{ISP}(a+1,b)$ with ground set of size $m+|V|$.
By Theorem~\ref{thm:fixed},
\[
 m+|V|\le n(a+1,b)
 =(1+o(1))\max_{0\le s\le a}
 (s+1)\binom{a+b-s}{b}.
\]
\end{proof}

For $a=b$, both constructions (Construction~\ref{constr1} and Construction~\ref{construct_ground} with $s=1$) asymptotically match the bound from the corollary, but distribute the sum differently. Construction~\ref{constr1} has $m=\binom{2a}{a},|V|=2a$, whereas Construction~\ref{construct_ground} with $s=1$ has $m=\frac12\binom{2a}{a}, |V|\sim \frac12\binom{2a}{a}$. This corollary provides weak stability: a system with almost the maximum number of pairs must have a relatively small ground set.

Another question, also considered by Tuza, is to look at the ground set of only the larger side of an ISP. Write $n_1(a,b)$ for this maximum $|\bigcup_i A_i|$ over all ISPs$(a,b)$. For $s=\lceil a/(b+1)\rceil$ Construction~\ref{construct_ground} satisfies $V=\bigcup_i A_i$, giving the following corollary.

\begin{corollary}\label{cor:one-sided}
For integers $a\ge b\ge1$ and $a\to\infty$, we have
\[
 n_1(a,b)\sim n(a,b)\sim K(a,b).
\]
\end{corollary}

One of the broader contributions of this paper is that we introduce a set of tools that allow for  a better understanding of the structure of large ISPs.

\subsection{Critical families and maximal intersecting families}

One of the motivations for studying ISP is the study of $\tau$-critical hypergraphs. A \emph{transversal} of a family $\F$ is a set meeting every member of $\F$. The minimum size of a transversal is denoted by $\tau(\F)$. We call $\F$ \emph{$\tau$-critical} if removing any set decreases its transversal number. Let $v_{\max}(k,t)$ be the maximum size of the ground set of a $\tau$-critical family $\F$ with $\tau(\F)=t$ and $|F|\le k$ for every $F\in\F$. The problem of determining this quantity goes back to Erd\H{o}s and Gallai~\cite{erdos-gallai}, who
proved that $v_{\max}(2,t)=2t$.

For any $k,t \ge 2$ Gy\'arf\'as, Lehel and Tuza~\cite{glt} proved the lower bound
\begin{equation}\label{eq:critical-lower}
 v_{\max}(k,t)\ge
 \binom{k+t-2}{k-1}+k+t-2.
\end{equation}
Their construction is the family $\F$ of all sets $A_i$ in
Construction~\ref{construct_ground} with $a=k$, $b=t-1$ and $s=1$.

For every fixed $k\ge2$, Gy\'arf\'as, Lehel and
Tuza~\cite{glt} proved that
$v_{\max}(k,t)=\Theta_k(t^{k-1})$ as $t\to\infty$.
For $k=3$, K\'ezdy and Lehel~\cite{kezdy-lehel-asymptotic}
established
$v_{\max}(3,t)=(1/2+o(1))t^2$, resolving asymptotically a conjecture of Szemer\'edi and Petruska~\cite{szemeredi-petruska}. Tuza~\cite{tuza-critical} asked whether equality holds in~\eqref{eq:critical-lower} for $t\ge k$, or at least whether this lower bound is asymptotically sharp.

The upper bound for $v_{\max}(k,t)$ is closely related to the quantity $n(a,b)$, through the following observation used by Tuza~\cite{tuza-critical}. If $\F$ is $\tau$-critical with $\tau(\F)=t$, then for each $F\in\F$ there is a $(t-1)$-set $T_F$ meeting every set in $\F\setminus\{F\}$. $T_F$ must avoid $F$, since otherwise it would be a transversal of $\F$. Thus the pairs $(F,T_F)$ form an ISP, and
\[
 v_{\max}(k,t)\le n_1(k,t-1).
\]

Together with Construction~\ref{construct_ground}, Theorem~\ref{thm:fixed} therefore gives the following asymptotic result.

\begin{corollary}\label{cor:critical}
For $2\le t\le k+1$, as $k\to\infty$,
\[
 v_{\max}(k,t)=(1+o(1))K(k,t-1).
\]
\end{corollary}

\begin{proof}
The upper bound follows from
$v_{\max}(k,t)\le n_1(k,t-1)\le n(k,t-1)$ and
Theorem~\ref{thm:fixed}.

For the lower bound, take $a=k$, $b=t-1$ and $s=\lceil k/t\rceil$ in Construction~\ref{construct_ground}, and let $\F$ consist of all $A_i$. Every set in $\binom{X}{t}$ is a transversal of $\F$, and no smaller transversal exists, so $\tau(\F)=t$. Deleting any member $S\cup P_S$ from $\F$ makes $X\setminus S$ a transversal of size $t-1$. Thus $\F$ is $\tau$-critical and has size of the ground set at least $K(k,t-1)$.
\end{proof}

Another closely related topic, initiated by Erd\H{o}s and Lov\'asz~\cite{erdos-lovasz}, is the study of maximal intersecting families. A $k$-uniform family $\F$ is \emph{intersecting} if
$F\cap G\ne\varnothing$ for all $F,G\in\F$. Such a family is called \emph{maximal} if, for any $k$-element set $F\notin\F$, the family $\F\cup\{F\}$ is not intersecting. Note that the ground set is not fixed here, so we allow $F$ to contain elements outside the ground set of $\F$. Every set from $\F$ is a transversal, so $\tau(\F)\le k$. If a transversal of size $k-1$ existed, we could extend it to a $k$-set using a new element and add this set to $\F$,
contradicting maximality. Thus $\tau(\F)=k$. Moreover, $\F$ coincides with the family of all its $k$-element transversals.

Erd\H{o}s and Lov\'asz asked how large $|\F|$ and the ground set of such a family can be. For the maximum number of sets, Frankl, Wang, Yang and Zhang~\cite{frankl-wang-yang-zhang}
recently improved the previous best lower bound of Frankl, Ota and Tokushige~\cite{frankl-ota-tokushige}.
Together with Zakharov's upper bound~\cite{zakharov-maximal}, this gives
\[
 \max_{\F}|\F|=k^k\exp\{-k^{1/2+o(1)}\}
 \qquad(k\to\infty).
\]

Let $N(k)$ be the maximum size of the ground set of a maximal intersecting $k$-uniform family.
Erd\H{o}s and Lov\'asz~\cite{erdos-lovasz} constructed a family showing that
\[
 N(k)\ge 2k-2+\frac12\binom{2k-2}{k-1}.
\]
Their construction is obtained from Construction~\ref{construct_ground} with $a=k$, $b=k-1$
and $s=1$ by using the same new element for $S$ and
$X\setminus S$. The family consists of the resulting sets $A_i$ together with all sets from $\binom{X}{k}$.

Tuza~\cite{tuza-critical} proved that
\[
 N(k)\le (2+o(1))\binom{2k-2}{k-1}.
\]
Majumder~\cite{majumder} showed that $N(k)\le \left(\frac32+o(1)\right)\binom{2k-2}{k-1}.$

As observed by Tuza~\cite{tuza-critical}, $N(k)\le n(k,k-1)$. Thus Theorem~\ref{thm:fixed} already gives $ N(k)\le(1+o(1))\binom{2k-2}{k-1},$ leaving an asymptotic factor of two between the lower and upper bounds. By using the additional structure of maximal intersecting families, we asymptotically close this gap. Our last main result shows that the construction of Erd\H{o}s and Lov\'asz is asymptotically optimal.

\begin{theorem}\label{thm:maximal}
As $k\to\infty$,
\[
 N(k)=\left(\frac12+o(1)\right)\binom{2k-2}{k-1}.
\]
\end{theorem}

The rest of the paper is organized as follows. Sections~\ref{sec:prelim} and~\ref{section_sampling} contain auxiliary
results. We prove Theorem~\ref{thm:fixed} in Section~\ref{sec:fixed}. In Section~\ref{sec:maximal}, we give a sketch and then a proof of
Theorem~\ref{thm:maximal}, assuming a technical lemma. We prove this technical lemma in Section~\ref{sec:singleton}.

\section{Preliminaries}\label{sec:prelim}

Families $\A,\B$ are \emph{cross-intersecting} if $A\cap B\ne\varnothing$ for every $A\in\A$, $B\in\B$.

Fix an ISP $\{(A_i,B_i):i\in[m]\}$ with ground set $V$.
For each $v\in V$, define
\[
 d_A(v)=|\{i:v\in A_i\}|,\qquad
 d_B(v)=|\{i:v\in B_i\}|,\qquad
 d(v)=d_A(v)+d_B(v).
\]
Since $A_i\cap B_i=\varnothing$, $d(v)$ is the number of pairs containing $v$. Double counting gives
\begin{equation}\label{eq:incidence}
 |V| = \sum_{i=1}^m\sum_{v \in (A_i \cup B_i)}\frac1{d(v)}.
\end{equation}

\subsection{A degree-weighted Tuza's inequality}

A family $\{(A_i,B_i):i\in[m]\}$ of pairs of finite sets is a \emph{weak ISP} if $A_i\cap B_i=\varnothing$ for every $i$, and at least one of $A_i\cap B_j$ and $A_j\cap B_i$ is nonempty whenever $i\ne j$. We first state Tuza's inequality for weak ISP.

\begin{theorem}[Tuza~\cite{tuza-weak}]\label{tuza-weak}
Let $0 < p<1$ be any real number and $q=1-p$. Then, for every weak ISP $\{(A_i,B_i):1\le i\le m\}$, we have
\begin{equation}\label{eq:tuza-colour}
 \sum_{i=1}^m p^{|A_i|}q^{|B_i|}\le1.
\end{equation}
\end{theorem}

For ISP, where both cross-intersections are required, we prove a degree-weighted refinement of Theorem~\ref{tuza-weak}. For $0<p<1$, put
\[
 M(p)=\max_{s\in\mathbb Z_{\ge1}}sp^s.
\]

\begin{lemma}\label{lem:colour}
Let $1/2 \le p< 1$ be any real number and $q=1-p$. Every ISP
$\{(A_i,B_i):1\le i\le m\}$ satisfies
\begin{equation}\label{eq:colour}
 \sum_{i=1}^m p^{|A_i|}q^{|B_i|}
       \sum_{v\in A_i\cup B_i}\frac1{d(v)}\le M(p).
\end{equation}
\end{lemma}

\begin{proof}
Choose a uniformly random ordering $\sigma$ of the elements in $[m]$. Call a pair $(v,A_i)$ \emph{free} if $v \in A_i$ and $i$ precedes every index $j$ with $v\in B_j$ in our random order $\sigma$. Similarly, call a pair $(v,B_i)$ free if $v \in B_i$ and $i$ follows every index $j$ with $v\in A_j$. Let $X_i$ and $Y_i$ count the free pairs with $A_i$ and $B_i$ respectively.

Independently of $\sigma$, choose a random set $S\subseteq V$
by including each element independently with probability $p$.
For each $i\in[m]$, let $E_i$ be the event that
$v\in S$ for every $v\in A_i$ such that $(v,A_i)$ is not free, and $v\notin S$ for every $v\in B_i$ such that $(v,B_i)$ is not free.

Conditional on the ordering,
\[
 \mathbb{P}(E_i\mid \sigma) =p^{|A_i|-X_i}q^{|B_i|-Y_i}.
\]

Events $E_i$ are pairwise disjoint. Indeed, if $i$ precedes $j$, choose $v\in A_j\cap B_i$. The pair $(v,A_j)$ is not free
because $i$ is earlier, and the pair $(v,B_i)$ is not free because $j$ is later. The events for $i$ and $j$ therefore make
opposite requirements on $v$: $v\in S$ and $v\notin S$. Summing and averaging gives
\begin{equation}\label{eq:free-events}
 \sum_i p^{|A_i|}q^{|B_i|}\EE[p^{-X_i}q^{-Y_i}]\le1.
\end{equation}

The probabilities of being first or last in the relevant order give
\[
 \EE X_i=\sum_{v\in A_i}\frac1{d_B(v)+1},\qquad
 \EE Y_i=\sum_{v\in B_i}\frac1{d_A(v)+1}.
\]
For $v\in A_i$ we have $d_B(v)+1\le d(v)$, and the analogous
inequality holds for $v\in B_i$. Also, for nonnegative integers $x,y$,
\[
 p^{-x}q^{-y}\ge p^{-(x+y)}\ge\frac{x+y}{M(p)}.
\]
Substitution in~\eqref{eq:free-events} gives 
$$M(p) \ge \sum_i p^{|A_i|}q^{|B_i|}\EE[X_i+Y_i]\ge \sum_i p^{|A_i|}q^{|B_i|}\Big(\sum_{v\in A_i}\frac1{d(v)}+\sum_{v\in B_i}\frac1{d(v)}\Big)$$
and proves~\eqref{eq:colour}.
\end{proof}

\subsection{Conflicting pairs in a random permutation}

Let $(A_1,B_1)$ and $(A_2,B_2)$ be pairs of subsets of a set $V$,
with $A_i\cap B_i=\varnothing$ for $i=1,2$.
Choose a uniformly random permutation of $V$, and let
$\widetilde E_i$ be the event that every element of $A_i$
precedes every element of $B_i$. Clearly, we have
\[
 \PP(\widetilde E_i)=\binom{|A_i|+|B_i|}{|A_i|}^{-1}.
\]

If both $A_1\cap B_2$ and $A_2\cap B_1$ are nonempty,
the events $\widetilde E_1$ and $\widetilde E_2$ are disjoint.
Indeed, for $x\in A_1\cap B_2$ and $y\in A_2\cap B_1$,
the event $\widetilde E_1$ requires $x$ to precede $y$,
whereas $\widetilde E_2$ requires the reverse.
Call the two pairs \emph{conflicting} if exactly one of
$A_1\cap B_2$ and $A_2\cap B_1$ is nonempty. We show below that, for conflicting pairs, the events
$\widetilde E_1$ and $\widetilde E_2$ rarely occur simultaneously.

\begin{lemma}\label{lem:conflict}
If the pairs $(A_1,B_1)$ and $(A_2,B_2)$ are conflicting, then
\begin{equation}\label{eq:conflict-bound}
 \PP(\widetilde E_2\mid\widetilde E_1)
 \le\frac1{\min\{|A_2|,|B_2|\}+1},
 \qquad
 \PP(\widetilde E_1\mid\widetilde E_2)
 \le\frac1{\min\{|A_1|,|B_1|\}+1}.
\end{equation}
\end{lemma}

\begin{proof}
Without loss of generality, assume that
$A_1\cap B_2=\varnothing$ and $A_2\cap B_1\ne\varnothing$.
Fix $v\in A_2\cap B_1$.
In every permutation belonging to
$\widetilde E_1\cap\widetilde E_2$, all elements of $A_1$
precede $v$, and $v$ precedes all elements of $B_2$.

For each such permutation and each $w\in\{v\}\cup B_2$,
swap the positions of $v$ and $w$.  The resulting permutation still belongs to $\widetilde E_1$:  no element of $A_1$ moves, and both swapped positions lie  after every element of $A_1$.  Moreover, $w$ is now the first element of $\{v\}\cup B_2$. Thus each permutation in $\widetilde E_1\cap\widetilde E_2$ gives $|B_2|+1$ distinct permutations in $\widetilde E_1$. No permutation is counted twice, so we have 
\[
 (|B_2|+1)\PP(\widetilde E_1\cap\widetilde E_2)
 \le\PP(\widetilde E_1),
\]
and hence
\[
 \PP(\widetilde E_2\mid\widetilde E_1)
 \le\frac1{|B_2|+1}.
\]

Similarly, swapping $v$ with any $w\in\{v\}\cup A_1$
preserves $\widetilde E_2$: no element of $B_2$ moves,
and both swapped positions lie before every element of $B_2$.
Therefore
\[
 (|A_1|+1)\PP(\widetilde E_1\cap\widetilde E_2)
 \le\PP(\widetilde E_2),
 \qquad
 \PP(\widetilde E_1\mid\widetilde E_2)
 \le\frac1{|A_1|+1}.
\]
These two bounds imply~\eqref{eq:conflict-bound}.
\end{proof}

\subsection{Sketch of the proof of Theorem~\ref{thm:fixed}}

We distinguish between elements of small and large degree.
For an ISP with ground set $V$, fix a constant $D$ and put
\[
 L=\{v\in V:d(v)\le D\},\qquad H=V\setminus L.
\]
We estimate $|L|$ and $|H|$ separately.

\begin{enumerate}
 \item Given an ISP$(a,b)$, sample each pair independently with probability $p$. Then delete from the ground set the elements of degree one in the sampled system, and apply Bollob\'as's inequality to the new ISP. This random restriction argument allows to bound $|H|$. The same argument allows us to discard pairs containing too many elements of $L$.

 \item Choose a uniformly random ordering of $H$, and let
 \[
 E_i=\{\text{all elements of }A_i\cap H
       \text{ precede all elements of }B_i\cap H\}.
 \]
If $(A_i\cap H,B_i\cap H)$ and $(A_j\cap H,B_j\cap H)$ are conflicting, Lemma~\ref{lem:conflict} shows that $\PP(E_j\mid E_i)$
is small. Among the indices for which $E_i$ occurs, remove those belonging to conflicting pairs. All cross-intersections between the remaining pairs lie in $L$.

 \item The selected pairs $(A_i\cap L,B_i\cap L)$ form an ISP, so we can apply Lemma~\ref{lem:colour}. It gives the required bound on $|L|$. We remark that  cleaning the conflicting pairs is crucial, since afterwards we get a product configuration: for each two pairs $(A_1,B_1), (A_2,B_2)$ from our ISP, either their restriction to $H$ or their restriction to $L$ gives a proper ISP pair.


\end{enumerate}

\section{Sampling and degree bounds}\label{section_sampling}

\subsection{Random sampling of pairs}

Throughout this section, fix an ISP $\{(A_i,B_i):i\in[m]\}$ with ground set $V$.
We call an element $v\in V$ \emph{private} if $d(v)=1$.
Note that removing private elements from ISP preserves intersection conditions.

\begin{lemma}\label{lem:sampling}
Let $0 < p<1$ be any real number and $q=1-p$. Every ISP$(a,b)$, satisfies
\begin{equation}\label{eq:side-sampling}
 \sum_{i=1}^{m}\exp\left(
   \frac12\log\frac{a+b}{a}\sum_{v\in A_i}q^{d(v)-1}
   +\frac12\log\frac{a+b}{b}\sum_{v\in B_i}q^{d(v)-1}
       \right)
 \le\frac{\binom{a+b}{b}}{p}.
\end{equation}
\end{lemma}
\begin{proof}
Sample each pair independently with probability $p$, and let
$I\subseteq[m]$ be the set of sampled indices.
For $i\in I$, let $X_i,Y_i$ count the private elements in
$A_i,B_i$, respectively, in the sampled system.
Delete all these elements, obtaining pairs $(A_i',B_i')$ with
\[
 |A_i'|=a-X_i,\qquad |B_i'|=b-Y_i.
\]
The pairs $(A_i',B_i')$ still form an ISP. By Bollob\'as's inequality,
\begin{equation}\label{eq:bollobas_short}
     1\ge\sum_{i\in I}\binom{|A_i'|+|B_i'|}{|B_i'|}^{-1}
 =\binom{a+b}{b}^{-1}
   \sum_{i\in I}
   \frac{\binom{a+b}{b}}
        {\binom{a+b-X_i-Y_i}{b-Y_i}}.    
\end{equation}

We have
\[
\begin{aligned}
 \frac{\binom{a+b}{b}}
      {\binom{a+b-X_i-Y_i}{b-Y_i}}
 &\ge\frac{\binom{a+b}{b}}{\binom{a+b-X_i}{b}}\\
 &=\prod_{j=0}^{X_i-1}\frac{a+b-j}{a-j}
 \ge\left(\frac{a+b}{a}\right)^{X_i}.
\end{aligned}
\]
Similarly, we obtain the lower bound $\bigl((a+b)/b\bigr)^{Y_i}$. Taking the geometric mean of these two bounds, we obtain 
\[
 \frac{\binom{a+b}{b}}
      {\binom{a+b-X_i-Y_i}{b-Y_i}}
 \ge\exp\left(
       \frac{X_i}{2}\log\frac{a+b}{a}
       +\frac{Y_i}{2}\log\frac{a+b}{b}
     \right).
\]
Substituting into inequality~\eqref{eq:bollobas_short} gives
\[
 \binom{a+b}{b}\ge\sum_{i\in I}
 \exp\left(
   \frac{X_i}{2}\log\frac{a+b}{a}
   +\frac{Y_i}{2}\log\frac{a+b}{b}
 \right).
\]

Conditional on $i\in I$, an element $v\in A_i\cup B_i$ is
private precisely when all other $d(v)-1$ pairs containing it
are not sampled. Thus
\[
 \EE[X_i\mid i\in I]=\sum_{v\in A_i}q^{d(v)-1},
 \qquad
 \EE[Y_i\mid i\in I]=\sum_{v\in B_i}q^{d(v)-1}.
\]
Taking expectations and applying Jensen's inequality, we obtain
\[
\begin{aligned}
 \binom{a+b}{b}
 &\ge\EE\sum_{i\in I}
 \exp\left(
   \frac{X_i}{2}\log\frac{a+b}{a}
   +\frac{Y_i}{2}\log\frac{a+b}{b}
 \right)\\
 &=p\sum_{i=1}^m
 \EE\left[
 \exp\left(
   \frac{X_i}{2}\log\frac{a+b}{a}
   +\frac{Y_i}{2}\log\frac{a+b}{b}
 \right)\,\middle|\,i\in I
 \right]\\
 &\ge p\sum_{i=1}^m
 \exp\left(
   \frac12\log\frac{a+b}{a}\,\EE[X_i\mid i\in I]
   +\frac12\log\frac{a+b}{b}\,\EE[Y_i\mid i\in I]
 \right)\\
 &=p\sum_{i=1}^m
 \exp\left(
   \frac12\log\frac{a+b}{a}\sum_{v\in A_i}q^{d(v)-1}
   +\frac12\log\frac{a+b}{b}\sum_{v\in B_i}q^{d(v)-1}
 \right).
\end{aligned}
\]
Dividing by $p$ proves~\eqref{eq:side-sampling}.
\end{proof}

Fix a constant $D$ and put
\[
 L=\{v\in V:d(v)\le D\},\qquad H=V\setminus L,\qquad
 s_i=|(A_i\cup B_i)\cap L|.
\]

The following corollary allows us to control pairs containing too many elements of small degree when $a$ and $b$ have the same order of magnitude.

\begin{corollary}\label{cor:trace-tail}
For every $\varepsilon>0$ and $0<\eta\le1/2$, there is an integer $T=T(D,\varepsilon,\eta)$ with the
following property. Let $\{(A_i,B_i):i\in[m]\}$ be an ISP$(a,b)$, where $a\ge b$, $a\ge2$ and $b/(a+b)\ge\eta$. Then
\begin{equation}\label{eq:trace-tail}
 \sum_{\substack{i\in[m]\\ s_i>T}} s_i \le\varepsilon\binom{a+b}{b}.
\end{equation}
\end{corollary}

\begin{proof}

Since $a\ge b$ and $b/(a+b)\ge\eta$,
\[
 \frac12\log\frac{a+b}{b}
 \ge\frac12\log\frac{a+b}{a}
 \ge\frac{b}{2(a+b)}\ge\frac\eta2.
\]
Applying Lemma~\ref{lem:sampling} with $p=q=1/2$, we obtain

\[
 2\binom{a+b}{b}
 \ge\sum_{i=1}^m
 \exp\left(\frac{\eta}{2}
       \sum_{v\in A_i\cup B_i}2^{1-d(v)}\right)
 \ge\sum_{i=1}^m
 \exp\left(\frac{\eta}{2^D}
       s_i \right).
\]

Choose $T=T(D, \varepsilon,\eta)$ large enough that
\[
 t\le\frac{\varepsilon}{2}
       \exp\left(\frac{\eta t}{2^D}\right)
 \qquad\text{for every integer }t>T.
\]
Consequently,
\[
 \varepsilon\binom{a+b}{b}
 \ge\frac{\varepsilon}{2}
   \sum_{i=1}^m
   \exp\left(\frac{\eta}{2^D}
       s_i \right) \ge  \sum_{\substack{i\in[m]\\ s_i>T}} s_i.
\]
\end{proof}

The next lemma also follows from Lemma~\ref{lem:sampling} and gives an upper bound on the number of elements of large degree.

\begin{lemma}\label{lem:high-degree}
For every $\varepsilon>0$, there is an integer $D=D(\varepsilon)$ such that every ISP$(a,b)$, where $a\ge b\ge1$ and $a\ge2$, satisfies
\begin{equation}\label{eq:high-degree}
 |H|
 \le\varepsilon\frac{a+b}{b}\binom{a+b}{b}.
\end{equation}
\end{lemma}

\begin{proof}
Since $a\ge b$ and $\log(1+x)\ge x/(1+x)$ for $x\ge0$, we have
$\frac12\log\frac{a+b}{b}\ge\frac12\log\frac{a+b}{a}
\ge\frac{b}{2(a+b)}$. Hence Lemma~\ref{lem:sampling}
and Jensen's inequality give, for every $0<p<1$,
\[
\begin{aligned}
 \frac{\binom{a+b}{b}}{p}
 &\ge\sum_{i=1}^m
 \exp\left(\frac{b}{2(a+b)}
       \sum_{v\in A_i\cup B_i}(1-p)^{d(v)-1}\right)\\
 &\ge m\exp\left(\frac{b}{2(a+b)m}
       \sum_{i=1}^m\sum_{v\in A_i\cup B_i}
       (1-p)^{d(v)-1}\right)\\
 &=m\exp\left(\frac{b}{2(a+b)m}
       \sum_{v\in V}d(v)(1-p)^{d(v)-1}\right).
\end{aligned}
\]
Taking logarithms and multiplying by $m$, we obtain
\begin{equation}\label{eq:log-sampling}
 m\log\frac{\binom{a+b}{b}}{pm}
 \ge\frac{b}{2(a+b)}
       \sum_{v\in V}d(v)(1-p)^{d(v)-1}.
\end{equation}

We integrate~\eqref{eq:log-sampling} over $0<p<1/D$.
For every $v\in H$,
\[
 \int_0^{1/D}d(v)(1-p)^{d(v)-1}\,dp
 =1-(1-1/D)^{d(v)}
 \ge1-e^{-d(v)/D}
 \ge1-e^{-1}.
\]
Since $d(v)(1-p)^{d(v)-1}\ge0$ for every $v\in L$,
we obtain
\[
 m\int_0^{1/D}\log\frac{\binom{a+b}{b}}{pm}\,dp
 \ge\frac{b(1-e^{-1})}{2(a+b)}\,|H|.
\]

On the other hand, Bollob\'as's inequality and
$\log x\le x-1$ give
\[
 m\le\binom{a+b}{b},
 \qquad
 m\log\frac{\binom{a+b}{b}}m
 \le\binom{a+b}{b}-m.
\]
Consequently,
\[
\begin{aligned}
 m\int_0^{1/D}\log\frac{\binom{a+b}{b}}{pm}\,dp
 &=m\int_0^{1/D}
   \left(\log\frac{\binom{a+b}{b}}m-\log p\right)\,dp\\
 &=\frac mD
   \left(\log\frac{\binom{a+b}{b}}m+1+\log D\right)\\
 &\le\frac{\binom{a+b}{b}+m\log D}{D}\\
 &\le\frac{1+\log D}{D}\binom{a+b}{b}.
\end{aligned}
\]
Combining these bounds yields
\[
 \frac{b(1-e^{-1})}{2(a+b)}\,|H|
 \le \frac{1+\log D}{D}\binom{a+b}{b}.
\]
Choose $D=D(\varepsilon)$ large enough that $\frac{2(1+\log D)}{(1-e^{-1})D}\le\varepsilon$. This gives the desired bound.
\end{proof}

\section{Proof of Theorem~\ref{thm:fixed}}\label{sec:fixed}

We first state the elementary estimates for $K(a,b)$.

\begin{lemma}\label{lem:optimization}
Let $a\ge b\ge1$ be integers, and put $p=a/(a+b)$.
There is an absolute constant $c>0$ such that
\begin{equation}\label{eq:K-comparison}
 c\frac ab\binom{a+b}{b}
 \le K(a,b)\le M(p)\binom{a+b}{b}.
\end{equation}
Moreover, for every fixed $\eta>0$, when
$b/(a+b)\ge\eta$, we have
\begin{equation}\label{eq:K-compact}
 K(a,b)=\bigl(M(p)+o(1)\bigr)\binom{a+b}{b}
 \qquad(a\to\infty).
\end{equation}
\end{lemma}

\begin{proof}
Recall that $K(a,b)=\max_{1\le s\le a}s\binom{a+b-s}{b}$.
Write $(u)_j=u(u-1)\cdots(u-j+1)$. Then
\begin{equation}\label{eq:binomial-product}
 \frac{K(a,b)}{\binom{a+b}{b}}
 =\max_{1\le s\le a}\frac{s(a)_s}{(a+b)_s}
 \le\max_{1\le s\le a}sp^s
 \le M(p).
\end{equation}
This proves the upper bound in~\eqref{eq:K-comparison}.

Suppose now that $b/(a+b)\ge\eta$, and choose an integer
$s\ge1$ such that $sp^s=M(p)$. By maximality,
\[
 sp^s\ge(s-1)p^{s-1},
\]
so $s(1-p)\le1$ and hence $s\le1/\eta$. Thus $s$ is bounded by a constant. Therefore,
\[
\begin{aligned}
 \frac{K(a,b)}{\binom{a+b}{b}}
 &\ge\frac{s(a)_s}{(a+b)_s}\\
 &=M(p)\prod_{j=0}^{s-1}
       \frac{1-j/a}{1-j/(a+b)}\\
 &=M(p)\bigl(1+o(1)\bigr).
\end{aligned}
\]
This and \eqref{eq:binomial-product} prove~\eqref{eq:K-compact}.

For the lower bound in~\eqref{eq:K-comparison}, first suppose
$a<2b$. Taking $s=1$ gives
\[
 K(a,b)\ge p\binom{a+b}{b}
 \ge\frac12\binom{a+b}{b}
 \ge\frac{a}{4b}\binom{a+b}{b}.
\]
If $a\ge2b$, take $s=\lfloor a/(2b)\rfloor$.
Then $s\ge a/(4b)$, and
\[
 \frac{\binom{a+b-s}{b}}{\binom{a+b}{b}}
 =\prod_{j=0}^{b-1}\left(1-\frac{s}{a+b-j}\right)
 \ge\left(1-\frac1{2b}\right)^b
 \ge\frac12.
\]
Consequently,
\[
 K(a,b)\ge s\binom{a+b-s}{b}
 \ge\frac{a}{8b}\binom{a+b}{b}.
\]
This proves the lower bound.
\end{proof}

\begin{proof}[Proof of Theorem~\ref{thm:fixed}]
Construction~\ref{construct_ground} gives $n(a,b)\ge K(a,b)$. We prove the upper bound. Fix $\varepsilon>0$, and let $\{(A_i,B_i):i\in[m]\}$ be an ISP with $|A_i|=a$ and $|B_i|=b$ for each $i \in [m]$.

For an integer $D\ge2$, put  $L=\{v\in V:d(v)\le D\}$ and $H=V\setminus L$. By Lemma~\ref{lem:high-degree} and the lower bound in~\eqref{eq:K-comparison}, we can choose $D=D(\varepsilon)$ large enough that
\begin{equation}\label{eq:fixed-high}
 |H|\le\varepsilon K(a,b).
\end{equation}
It remains to bound $|L|$.

For each $i\in[m]$, write
\[
 x_i=|A_i\cap L|,\qquad
 y_i=|B_i\cap L|,\qquad
 s_i=x_i+y_i.
\]

Applying Lemma~\ref{lem:sampling} with $p=q=1/2$, using $2^{1-d(v)}\ge 2^{1-D}$ for $v\in L$ and
$\log\frac{a+b}{a} \ge \frac{b}{a+b}$, we obtain
\begin{equation}\label{eq:side-tail}
 2\binom{a+b}{b}
 \ge\sum_{i=1}^m
 \exp\left(c_1\left[
      \frac{b}{a+b}x_i
      +\log\frac{a+b}{b}\,y_i
 \right]\right),
\end{equation}
where $c_1>0$ depends only on $D$.

Since $\log((a+b)/b)\ge b/(a+b)$, for $y_i>0$ we have
\[
 \frac{b}{a+b}x_i+\log\frac{a+b}{b}\,y_i
 \ge\frac{b}{2(a+b)}s_i+\frac12\log\frac{a+b}{b}.
\]
Consequently,
\[
\begin{aligned}
 2\binom{a+b}{b}
 &\ge\sum_{i:y_i>0}
 \exp\left(c_1\left[
      \frac{b}{a+b}x_i
      +\log\frac{a+b}{b}\,y_i
 \right]\right)\\
 &\ge
 \left(\frac{a+b}{b}\right)^{c_1/2}
 \sum_{i:y_i>0}
 \exp\left(\frac{c_1b}{2(a+b)}s_i\right)\\
 &\ge
 \left(\frac{a+b}{b}\right)^{c_1/2}
 \frac{ec_1b}{2(a+b)}
 \sum_{i:y_i>0}s_i.
\end{aligned}
\]
Combining this with the lower bound
in~\eqref{eq:K-comparison}, we obtain
\begin{equation}\label{eq:bad-short-side}
 \sum_{i:y_i>0}s_i
 \le c_2\left(\frac{b}{a+b}\right)^{c_1/2}K(a,b),
\end{equation}
where $c_2>0$ also depends only on $D$.

Choose $0<\eta<1/2$ so that $c_2\eta^{c_1/2}\le\varepsilon$. This choice depends only on $\varepsilon$. We now consider two cases.

\medskip\noindent
\emph{Case 1: $b/(a+b)<\eta$.}

By~\eqref{eq:bad-short-side},
\[
 \sum_{i:y_i>0}s_i\le\varepsilon K(a,b).
\]
For the remaining indices, $B_i\cap L=\varnothing$ holds.
Thus deleting elements of $L$ from each $A_i$ preserves all cross-intersections among these pairs, and Bollob\'as's inequality gives
\[
 1\ge
 \sum_{i:y_i=0}\binom{a+b-x_i}{b}^{-1}.
\]
By the definition of $K(a,b)$,
\[
 x_i\binom{a+b-x_i}{b}\le K(a,b).
\]
Hence
\[
 K(a,b)
 \ge K(a,b)\sum_{i:y_i=0}\binom{a+b-x_i}{b}^{-1}
 \ge\sum_{i:y_i=0}x_i.
\]
Every element of $L$ occurs in at least one of the original
pairs, so
\[
 |L|
 \le\sum_{i:y_i=0}x_i+\sum_{i:y_i>0}s_i
 \le(1+\varepsilon)K(a,b).
\]
Together with~\eqref{eq:fixed-high}, this yields
\[
 |V|\le(1+2\varepsilon)K(a,b).
\]

\medskip\noindent
\emph{Case 2: $b/(a+b)\ge\eta$.}

In this case both $a$ and $b$ tend to infinity.
By Corollary~\ref{cor:trace-tail}, there exists an integer $T=T(D,\varepsilon,\eta)$ such that
\[
 \sum_{i:s_i>T}s_i\le\varepsilon\binom{a+b}{b}.
\]
Let $J=\{i\in[m]:s_i\le T\}$.

For each $i \in [m]$, put
\[
 \omega_i=
 \sum_{v\in(A_i\cup B_i)\cap L}\frac1{d(v)},
\]
where the degrees are taken in the original system.
Then
\[
 \sum_{i=1}^m\omega_i=|L|,
 \qquad
 \sum_{i\notin J}\omega_i
 \le\sum_{i\notin J}s_i
 \le\varepsilon\binom{a+b}{b}.
\]

For each $i\in J$, put
\[
 A_i'=A_i\cap H,\qquad B_i'=B_i\cap H.
\]
Choose a uniformly random permutation of $H$,
and let $E_i$ be the event that every element of $A_i'$
precedes every element of $B_i'$. Then
\[
 \PP(E_i)=\binom{a+b-x_i-y_i}{a-x_i}^{-1}.
\]
We call distinct indices $i,j\in J$ conflicting if the pairs
$(A_i',B_i')$ and $(A_j',B_j')$ are conflicting.

We first bound the number of indices conflicting with a fixed $i$.
For every such index $j$, one of $A_i'\cap B_j'$ and
$A_j'\cap B_i'$ is empty. The corresponding original
cross-intersection is nonempty and lies entirely in $L$.
Thus the pair $(A_j,B_j)$ contains an element of
$(A_i\cup B_i)\cap L$.
There are at most $T$ such elements, and each belongs to at most
$D-1$ other pairs. Hence $i$ conflicts with at most $(D-1)T$
indices.

Since $s_j\le T$ for every $j\in J$,
\[
 |A_j'|\ge a-T\ge b-T,\qquad |B_j'|\ge b-T.
\]
For a conflicting pair $i,j$ Lemma~\ref{lem:conflict} gives 
\[
 \PP(E_j\mid E_i)\le\frac1{b-T+1}=o(1).
\]

Let $I$ be the set of all indices $i\in J$ for which $E_i$ occurs but $E_j$ does not occur for any index $j\in J$ conflicting with $i$. By the union bound,
\begin{equation}\label{eq:fixed-select}
\begin{aligned}
 \PP(i\in I)
 &\ge\PP(E_i)\left(
       1-\sum_{\substack{j\in J\\ i,j\text{ conflicting}}}
       \PP(E_j\mid E_i)\right)\\
 &\ge
 \left(1-\frac{(D-1)T}{b-T+1}\right)
 \binom{a+b-x_i-y_i}{a-x_i}^{-1}\\
 &=(1-o(1))
 \binom{a+b-x_i-y_i}{a-x_i}^{-1}.
\end{aligned}
\end{equation}

Let $i,j\in I$ be distinct. By the definition of $I$, both $E_i$ and $E_j$ occur, so $A_i'\cap B_j'$ and $A_j'\cap B_i'$ cannot
both be nonempty. They cannot have exactly one nonempty intersection either, since we do not take conflicting pairs in $I$. Thus both intersections $A_i'\cap B_j'$ and $A_j'\cap B_i'$ are empty. All original cross-intersections between pairs from $I$ are consequently nonempty subsets of $L$, so
\[
 \{(A_i\cap L,B_i\cap L):i\in I\}
\]
is an ISP.

Apply Lemma~\ref{lem:colour} to this ISP with $p=\frac{a}{a+b}$. Taking expectations and using~\eqref{eq:fixed-select}, we obtain
\begin{equation}\label{eq:fixed-average}
\begin{aligned}
 M\left(\frac{a}{a+b}\right)
 &\ge
 \EE\sum_{i\in I}
 \left(\frac{a}{a+b}\right)^{x_i}
 \left(\frac{b}{a+b}\right)^{y_i}\omega_i\\
 &=
 \sum_{i\in J}\PP(i\in I)
 \left(\frac{a}{a+b}\right)^{x_i}
 \left(\frac{b}{a+b}\right)^{y_i}\omega_i\\
 &\ge
 (1-o(1))\sum_{i\in J}
 \frac{
       \left(\frac{a}{a+b}\right)^{x_i}
       \left(\frac{b}{a+b}\right)^{y_i}}
      {\binom{a+b-x_i-y_i}{a-x_i}}
 \,\omega_i.
\end{aligned}
\end{equation}

Recall that $x_i+y_i\le T$ and $a,b\ge\eta(a+b)$. Since $\eta$ and $T$ are fixed,
\begin{equation}\label{eq:fixed-cancellation}
\begin{aligned}
 \binom{a+b-x_i-y_i}{a-x_i}
 &=\binom{a+b}{b}
   \frac{(a)_{x_i}(b)_{y_i}}{(a+b)_{x_i+y_i}}\\
 &=(1+o(1))\binom{a+b}{b} \frac{a^{x_i}b^{y_i}}{(a+b)^{x_i+y_i}}.
\end{aligned}
\end{equation}
Substituting this into~\eqref{eq:fixed-average} and
using~\eqref{eq:K-compact}, we conclude that
\[
 \sum_{i\in J}\omega_i
 \le (1+o(1))
 M\left(\frac{a}{a+b}\right)\binom{a+b}{b}
 =(1+o(1))K(a,b).
\]

 For the next calculation, we need a lower bound on  $K(a,b)$.  Taking $s=1$ in the definition of $K(a,b)$ gives
\[
 K(a,b)\ge\binom{a+b-1}{b}
 =\frac{a}{a+b}\binom{a+b}{b}
 \ge\frac12\binom{a+b}{b}.
\]

Finally, we obtain
\[
\begin{aligned}
 |V|
 &=|H|
   +\sum_{i\in J}\omega_i
   +\sum_{i\notin J}\omega_i\\
 &\le
 \varepsilon K(a,b)
 +(1+o(1))K(a,b)
 +\varepsilon\binom{a+b}{b}\\
 &\le(1+3\varepsilon+o(1))K(a,b).
\end{aligned}
\]

Since $\varepsilon>0$ was arbitrary,
the theorem follows.
\end{proof}

\section{Proof of Theorem~\ref{thm:maximal}}\label{sec:maximal}

For cross-intersecting families $\A,\B\subseteq2^U$ let
\[
 R(\A,\B)=\{v\in U:A\cap B=\{v\}\text{ for some }A\in\A,B\in\B\}.
\]
We call such a pair $A,B$ a \emph{singleton witness} for $v$.
For each $v\in R(\A,\B)$, define
\begin{equation}\label{eq:W}
 w(v)=\max_{\substack{A\in\A,\ B\in\B\\A\cap B=\{v\}}}
          (|A|+|B|-2).
\end{equation}
The following theorem is a key technical ingredient in the proof.
\begin{theorem}\label{thm:singleton}
For any two cross-intersecting families $\A,\B\subseteq2^U$,
\begin{equation}\label{eq:singleton}
 \sum_{v\in R(\A,\B)}2^{-w(v)}\le1.
\end{equation}
\end{theorem}

We postpone the proof of Theorem~\ref{thm:singleton} until Section~\ref{sec:singleton}. We first prove the upper bound in Theorem~\ref{thm:maximal} using this inequality.

\subsection{Sketch of the proof}

Let $\F$ be a maximal intersecting $k$-uniform family.

\begin{enumerate}
 \item After reducing to a $\tau$-critical family, look at the $\operatorname{ISP}(k,k-1)$ formed by the sets of $\F$ and $(k-1)$-transversals. As in the proof of Theorem~\ref{thm:fixed}, remove elements of large degree and then elements in pairs containing too many elements of small degree. We then focus on the remaining elements, which form $L$,
and put $H=V\setminus L$.

\item For every $v\in L$, choose $A_v,B_v\in\F$ with
$A_v\cap B_v=\{v\}$, maximizing $|A_v\cap L|+|B_v\cap L|$. Consider both $(A_v,B_v)$ and $(B_v,A_v)$ and denote these pairs by $(A_i,B_i)$. Since $\F$ is intersecting, we have $A_i\cap B_j\ne\varnothing$ for all $i,j$.

\item Choose a uniformly random ordering of $H$, and let $E_i$ be the event that all elements of $A_i\cap H$ precede all elements of $B_i\cap H$. Among the indices for which $E_i$ occurs, remove both indices of every conflicting pair. For the remaining indices $I$, the families $\{A_i\cap L:i\in I\}$ and $\{B_i\cap L:i\in I\}$ are cross-intersecting.

\item Apply Theorem~\ref{thm:singleton} to
$\{A_i\cap L:i\in I\}$ and $\{B_i\cap L:i\in I\}$,
and take expectations to obtain the required bound on $|L|$.

\end{enumerate}

\subsection{Proof of the upper bound assuming Theorem~\ref{thm:singleton}}

It is enough to prove the upper bound for intersecting $k$-uniform families with $\tau = k$. Choose such a family $\F=\{F_1,\ldots,F_m\}$ with ground set
$V$ of maximum size and, among all such families, with $m$ as small as possible. Then no $F\in \mathcal F$  can have a private element $x$: otherwise, $F\setminus \{x\}$ is a transversal. In particular, removing any set from $\F$ does not change its ground set. Therefore, by minimality, $\F$ must be $\tau$-critical. 
For each $i\in[m]$, choose a $(k-1)$-transversal $T_i\subseteq V$ of $\F\setminus\{F_i\}$; since $\tau(\F)=k$,
we have $T_i\cap F_i=\varnothing$, and hence $\{(F_i,T_i):i\in[m]\}$ is an ISP$(k,k-1)$.

For $v\in V$, let $d(v)=|\{i:v\in F_i\}|$. Throughout the rest of the proof, we assume that $k$ is sufficiently large.

\begin{claim}\label{cl:regularize}
For every $\varepsilon>0$, there exist integers
$D=D(\varepsilon)$ and $T=T(\varepsilon)$ and a partition
$V=L\sqcup H$ such that $$|H|\le\varepsilon\binom{2k-2}{k-1}$$ and the following holds for $L$: $d(v)\le D$ for every $v\in L$ and $|F\cap L|\le T$ for every $F\in\F$.
\end{claim}

\begin{proof}
Apply Lemma~\ref{lem:high-degree} to the ISP
$\{(F_i,T_i):i\in[m]\}$ with $a=k$ and $b=k-1$. We can choose $D=D(\varepsilon)$ such that, for
\[
 L_0=\{v\in V:|\{i:v\in F_i\cup T_i\}|\le D\},
 \qquad H_0=V\setminus L_0,
\]
we have
\[
 |H_0|\le\frac{\varepsilon}{2}\binom{2k-2}{k-1}.
\]
In particular, $d(v)\le D$ for every $v\in L_0$.

Put $s_i=|(F_i\cup T_i)\cap L_0|$. By Corollary~\ref{cor:trace-tail}, there exists an integer $T=T(D,\varepsilon)$ such that
\[
 \sum_{i:s_i>T}s_i
 \le\frac{\varepsilon}{2}\binom{2k-2}{k-1}.
\]
Define
\[
 L=L_0\setminus\bigcup_{i:s_i>T}(F_i\cup T_i),
 \qquad H=V\setminus L.
\]
Then
\[
 |H|\le |H_0|+\sum_{i:s_i>T}s_i
 \le\varepsilon\binom{2k-2}{k-1}.
\]
It is easy to see that $L, H$ satisfy the conditions in the claim.
\end{proof}

We now focus on the elements of $L$ and keep track of the
intersections $F\cap L$ for $F\in\F$. For $v\in L$, define
\[
 w_L(v)=\max_{\substack{F,G\in\F\\F\cap G=\{v\}}}
          \bigl(|F\cap L|+|G\cap L|-2\bigr).
\]
For each $v\in L$, choose $A_v,B_v\in\F$ attaining this maximum.
Include both $(A_v,B_v)$ and $(B_v,A_v)$ in an indexed family
$\{(A_i,B_i):i\in [\ell]\}$, so that $\ell=2|L|$.

For $i\in [\ell]$, let $v_i$ be the unique element of $A_i\cap B_i$, and put
\[
 x_i=|A_i\cap L|,\qquad y_i=|B_i\cap L|.
\]
Thus
\begin{equation}\label{eq:chosen-witness}
 A_i\cap B_i=\{v_i\},\qquad
 x_i+y_i-2=w_L(v_i),\qquad 1\le x_i,y_i\le T.
\end{equation}
All cross-intersections $A_i\cap B_j$ are nonempty, since
$A_i,B_j\in\F$.

Put
\[
 A_i'=A_i\cap H,\qquad B_i'=B_i\cap H.
\]
Since $A_i$ and $B_i$ intersect only at $v_i\in L$, the sets
$A_i'$ and $B_i'$ are disjoint. We have $|A_i'| = k-x_i$ and
$|B_i'| = k-y_i$. Choose a uniformly random permutation of $H$, and let $E_i$ be the event that every element of $A_i'$ precedes
every element of $B_i'$. Then
\[
 \PP(E_i)=\binom{2k-x_i-y_i}{k-x_i}^{-1}.
\]

Call two distinct indices $i,j$ conflicting if exactly one of
$A_i'\cap B_j'$ and $A_j'\cap B_i'$ is nonempty.
Since $|A_j'|,|B_j'|\ge k-T$, Lemma~\ref{lem:conflict} gives
\[
 \PP(E_j\mid E_i)
 \le\frac{1}{\min\{|A_j'|,|B_j'|\}+1}
 \le\frac{1}{k-T+1} =o(1)
\]
whenever $i,j$ are conflicting.

We next bound the number of indices conflicting with a fixed $i$. A set $F\in\F$ can be chosen as $A_v$ or $B_v$  only for an element  $v\in F\cap L$. Taking into account that each pair $A_v,B_v$ appears in both orientations, the set $F$ therefore appears in at most $2|F\cap L|\le 2T$ indexed pairs. Each element of $L$ belongs to at most $D$ sets in $\F$, and hence to at most $2DT$ indexed pairs.

If $i,j$ are conflicting, then either 
$A_i\cap B_j$ or $A_j\cap B_i$ lies entirely in $L$.
The respective intersection is nonempty, and so $(A_j\cup B_j)$ contains an element of $(A_i\cup B_i)\cap L$. Since
\[
 |(A_i\cup B_i)\cap L|=x_i+y_i-1\le2T-1,
\]
there are at most $2DT(2T-1)$ indices conflicting with $i$.

Going back to our random experiment, take all indices $i$ for which $E_i$ occurs.
Remove every index belonging to a conflicting pair $i,j$ for which both $E_i$ and $E_j$ occur, and let $I \subset [\ell]$ be the set of remaining indices. By the union bound,
\begin{equation}\label{eq:mif-select}
\begin{aligned}
 \PP(i\in I)
 &\ge \PP(E_i)\left(
       1-\sum_{\substack{j\in [\ell]\\i,j\text{ conflicting}}}
              \PP(E_j\mid E_i)\right)\\
 &\ge \bigl(1-o(1)\bigr)
       \binom{2k-x_i-y_i}{k-x_i}^{-1}.
\end{aligned}
\end{equation}

Let $i,j\in I$ be distinct. Since both $E_i$ and $E_j$ occur,
the intersections $A_i'\cap B_j'$ and $A_j'\cap B_i'$ cannot
both be nonempty. If exactly one were nonempty, then $i$ and $j$
would be conflicting, contrary to the definition of $I$. Thus both
intersections are empty. Since $A_i'\cap B_i'=\varnothing$ as well,
we obtain
\begin{equation}\label{eq:mif-outside}
 (A_i\cap H)\cap(B_j\cap H)=\varnothing
 \qquad(i,j\in I).
\end{equation}

The events corresponding to $(A_v,B_v)$ and
$(B_v,A_v)$ are disjoint: the two nonempty sets
$A_v\cap H$ and $B_v\cap H$ cannot each precede
the other. Thus at most one of the two indices associated
with $v$ belongs to $I$, so for $i \in I$ the elements $v_i$ are distinct.

For our random permutation, consider the families
\[
 \A_I=\{A_i\cap L:i\in I\},\qquad
 \B_I=\{B_i\cap L:i\in I\}.
\]
By~\eqref{eq:mif-outside},
\[
 (A_i\cap L)\cap(B_j\cap L)=A_i\cap B_j\ne\varnothing
 \qquad(i,j\in I),
\]
so these families are cross-intersecting. Applying Theorem~\ref{thm:singleton} to $\A_I$ and $\B_I$, we obtain
\begin{equation}\label{eq:mif-trace-inequality}
 \sum_{i\in I}2^{-w_L(v_i)}\le1.
\end{equation}

Taking expectations in~\eqref{eq:mif-trace-inequality} and using~\eqref{eq:mif-select}, we obtain 
$$
 1\ge\EE\sum_{i\in I}2^{-w_L(v_i)}
   =\sum_{i\in [\ell]}\PP(i\in I)2^{-(x_i+y_i-2)} \ge\bigl(1-o(1)\bigr)
       \sum_{i\in [\ell]}
       \frac{2^{-(x_i+y_i-2)}}
            {\binom{2k-x_i-y_i}{k-x_i}}.
$$

Since $x_i,y_i \le T$, we have $\binom{2k-x_i-y_i}{k-x_i}\sim \binom{2k-2}{k-1} 2^{-(x_i+y_i-2)}$. Therefore, 
$$ \sum_{i=1}^{\ell}
 \frac{2^{-(x_i+y_i-2)}}{\binom{2k-x_i-y_i}{k-x_i}}
 \ge (1-o(1))\ell\binom{2k-2}{k-1}^{-1}
 = (1-o(1))\,2|L|\binom{2k-2}{k-1}^{-1}.
$$

We obtain
\[
 |L|\le\left(\frac12+o(1)\right)\binom{2k-2}{k-1}.
\]
Adding the bound on $|H|$
from Claim~\ref{cl:regularize}, we conclude that
\[
 |V|\le\left(\frac12+\varepsilon+o(1)\right)
          \binom{2k-2}{k-1}.
\]
Letting first $k\to\infty$ and then $\varepsilon\to0$
proves the desired upper bound.

\section{Proof of Theorem~\ref{thm:singleton}}\label{sec:singleton}

We induct on $|U|$. Set $2^{-w(v)}=0$ for $v\notin R(\A,\B)$.

For a set $K\subseteq U$, choose a uniformly random
subset $S\subseteq K$ and define
\[
\begin{aligned}
 \A'&=\{A\setminus K:A\in\A,\ A\cap K\subseteq S\},\\
 \B'&=\{B\setminus K:B\in\B,\ B\cap K\subseteq K\setminus S\}.
\end{aligned}
\]
If $A\cap K\subseteq S$ and $B\cap K\subseteq K\setminus S$,
then $A\cap B\subseteq U\setminus K$. Thus $\A',\B'$ are
cross-intersecting. For $v\in R(\A',\B')$, put
\[
 w'(v)=
 \max_{\substack{A'\in\A',\,B'\in\B'\\A'\cap B'=\{v\}}}
       (|A'|+|B'|-2),
\]
and set $2^{-w'(v)}=0$ for $v\notin R(\A',\B')$.

\begin{claim}\label{cl:singleton-restriction}
For every nonempty $K\subseteq U$,
\begin{equation}\label{eq:singleton-restriction}
 \sum_{v\in U}2^{-w(v)}
 \le
 \sum_{v\in K}2^{-w(v)}
 +\PP\bigl(R(\A',\B')\ne\varnothing\bigr).
\end{equation}
\end{claim}

\begin{proof}
Fix $v\in R(\A,\B)\setminus K$. Choose $A\in\A$ and $B\in\B$ with $A\cap B=\{v\}$ so that $t=|(A\cup B)\cap K|$ is as small as possible. Since $A\cap K$ and $B\cap K$
are disjoint,
\[
 \PP(A\cap K\subseteq S,\ B\cap K\subseteq K\setminus S)
 =2^{-t}.
\]

By the minimality of $t$, we have $w'(v)\le w(v)-t$ whenever $v\in R(\mathcal A',\mathcal B')$. Therefore,
\[
 \EE\,2^{-w'(v)}
 \ge 2^{-t}2^{-(w(v)-t)}
 =2^{-w(v)}.
\]
Summing over $v\in U\setminus K$ gives
\[
 \sum_{v\in U\setminus K}2^{-w(v)}
 \le \EE\sum_{v\in U\setminus K}2^{-w'(v)}.
\]
By induction, the sum inside the expectation is at most one,
and it is zero when $R(\A',\B')=\varnothing$. Consequently,
\[
 \sum_{v\in U\setminus K}2^{-w(v)}
 \le \PP\bigl(R(\A',\B')\ne\varnothing\bigr),
\]
which proves the claim.
\end{proof}

Without loss of generality, let $A_1\in\A$ be a set of minimum size among all sets in $\A\cup\B$, and put $a=|A_1|$.
Apply Claim~\ref{cl:singleton-restriction} for $K=A_1$.
If $S=A_1$, then $\B'=\varnothing$, since every $B\in\B$
has a nonempty intersection with $A_1$. Thus we have $\PP\bigl(R(\A',\B')\ne\varnothing\bigr)\le1-2^{-a},$ and
\begin{equation}\label{eq:singleton-smallest-set}
 \sum_{v\in U}2^{-w(v)}
 \le \sum_{v\in A_1}2^{-w(v)}+1-2^{-a}.
\end{equation}

By minimality of $|A_1|$, we have $w(v)\ge2a-2$ for every
$v\in R(\A,\B)$. If $a\ge4$, then
\[
 \sum_{v\in A_1}2^{-w(v)}
 \le \frac{a}{2^{2a-2}}
 \le \frac1{2^a},
\]
where the last inequality follows from $a\le2^{a-2}$.
Together with~\eqref{eq:singleton-smallest-set}, this proves
the theorem for $a\ge4$.

If $a=1$, write $A_1=\{v\}$. Every $B\in\B$ contains $v$,
so $v\in R(\A,\B)$. If $w(v)\ge1$, then
\eqref{eq:singleton-smallest-set} gives
\[
 \sum_{u\in U}2^{-w(u)}
 \le 2^{-w(v)}+\frac12\le1.
\]
If $w(v)=0$, then both $\A$ and $\B$ consist only of the
singleton $\{v\}$, and the required inequality is immediate.

It remains to consider $a\in\{2,3\}$. If $w(v)\ge2a-1$
for every $v\in R(\A,\B)$, then
\[
 \sum_{v\in A_1}2^{-w(v)}
 \le\frac{a}{2^{2a-1}}
 \le\frac1{2^a},
\]
so~\eqref{eq:singleton-smallest-set} again proves the result.
We may therefore assume that
\[
 H=\{v\in R(\A,\B):w(v)=2a-2\}
\]
is nonempty. Every singleton witness for an element of $H$
consists of two sets of size $a$.

Choose subfamilies $\A_0\subseteq\A$ and $\B_0\subseteq\B$
of sets of size $a$, minimal under deletion of sets subject
to having a singleton witness for every $v\in H$.
For each $A\in\A_0$, minimality supplies $v_A\in H$
such that every singleton witness for $v_A$ in
$(\A_0,\B_0)$ uses $A$. Choose $B_A\in\B_0$ with
$A\cap B_A=\{v_A\}$. Then
\[
 \{(A,B_A\setminus\{v_A\}):A\in\A_0\}
\]
is an $\operatorname{ISP}(a,a-1)$. The same construction applies with $B \in \B_0$. We will use the following estimate for these two ISPs.

\begin{claim}\label{cl:singleton-small-sets}
Let $a\in\{2,3\}$, and let
$\{(C_i,T_i):i\in[m]\}$ be an $\operatorname{ISP}(a,a-1)$.
For a uniformly random subset $S$ of its ground set,
\begin{equation}\label{eq:singleton-small-sets}
 \PP(C_i\subseteq S\text{ for some }i\in[m])
 \ge \frac{|\bigcup_{i=1}^m C_i|}{2^{2a-1}},
\end{equation}
unless $a=2$, $m=2$ and $C_1\cap C_2=\varnothing$.
\end{claim}

\begin{proof}
Suppose first that $a=2$, and write $T_i=\{t_i\}$. Since $|C_i|=2$, we have $m\le3$. If $m=1$, then
\[
 \PP(C_1\subseteq S)=\frac14=\frac{|C_1|}{8}.
\]
If $m=2$ and $C_1\cap C_2\ne\varnothing$, then
$|C_1\cap C_2|=1$, and
\[
 \PP(C_1\subseteq S\text{ or }C_2\subseteq S)
 =\frac14+\frac14-\frac18
 =\frac38
 =\frac{|C_1\cup C_2|}{8}.
\]
If $m=3$, then $C_i=\{t_1,t_2,t_3\}\setminus\{t_i\}$. Consequently,
\[
 \PP(C_i\subseteq S\text{ for some }i\in[3])
 =\PP(|S\cap\{t_1,t_2,t_3\}|\ge2)
 =\frac12\ge\frac38.
\]

Now suppose that $a=3$. We induct on $m$.
For $m=1$, the required bound is $1/8\ge3/32$.
Call an element \emph{unique} if it belongs to exactly one
of $C_1,\ldots,C_m$; it may still belong to some $T_j$.

If some $C_i$ has no unique elements, remove $(C_i,T_i)$.
Then $\bigcup_{j\ne i}C_j=\bigcup_{j=1}^m C_j$, and the required bound follows by induction.

Suppose that some $C_i$ has exactly one unique element. The event $C_i\subseteq S$, $T_i\cap S=\varnothing$ has probability $1/32$. On this event, no $C_j$ with $j\ne i$
is contained in $S$, since $C_j\cap T_i\ne\varnothing$.
Thus, by induction,
\[
\begin{aligned}
 \PP(C_j\subseteq S\text{ for some }j\in[m])
 &\ge
 \PP(C_j\subseteq S\text{ for some }j\ne i)+\frac1{32}\\
 &\ge
 \frac{|\bigcup_{j\ne i}C_j|+1}{32}
 =\frac{|\bigcup_{j=1}^m C_j|}{32}.
\end{aligned}
\]

We may now assume that every $C_i$ has at least two unique
elements. The sets $C_1,\ldots,C_m$ therefore split into groups such that sets from different groups are disjoint, and for each group of size at least two there is an element $v$ such that $C_i\cap C_j=\{v\}$ for any two distinct sets in that group. Let $g$ be the number of obtained groups. 

If $C_i$ belongs to a group of size $t$, then $T_i$ must contain $t-1$ distinct elements from the other sets in this group and at least one element from each other group. Thus $t+g-2\le|T_i|=2$. Every group has at most $4-g$ sets, giving
\[
 m\le g(4-g)\le4.
\]

For distinct $i,j$, we have $|C_i\cup C_j|\ge5$. Hence
\[
\begin{aligned}
 \PP(C_i\subseteq S\text{ for some }i\in[m])
 &\ge
 \sum_{i=1}^m\PP(C_i\subseteq S)
 -\sum_{1\le i<j\le m}\PP(C_i\cup C_j\subseteq S)\\
 &\ge\frac m8-\frac{\binom m2}{32}.
\end{aligned}
\]
If $m\le3$, then
\[
 \frac m8-\frac{\binom m2}{32}
 \ge\frac{3m}{32}
 \ge\frac{|\bigcup_{i=1}^m C_i|}{32}.
\]
If $m=4$, the inequality $m\le g(4-g)$ forces $g=2$,
with two sets in each group. The union of each group
has size five, so
\[
 \left|\bigcup_{i=1}^4 C_i\right|=10.
\]
The same estimate gives
\[
 \PP(C_i\subseteq S\text{ for some }i\in[4])
 \ge\frac48-\frac6{32}
 =\frac{10}{32}.
\]
This completes the proof.
\end{proof}

Let $K$ be the union of all sets in $\A_0\cup\B_0$. Since
$w(v)=2a-2$ for $v\in H$ and $w(v)\ge2a-1$ for
$v\in R(\A,\B)\setminus H$, we have
\begin{equation}\label{eq:singleton-weight-on-K}
 \sum_{v\in K}2^{-w(v)}
 \le\frac{|H|}{2^{2a-2}}
    +\frac{|K|-|H|}{2^{2a-1}}
 =\frac{|K|+|H|}{2^{2a-1}}.
\end{equation}

Choose a uniformly random subset $S\subseteq K$ and define
$\A',\B'$ as in Claim~\ref{cl:singleton-restriction}.
If $A\subseteq S$ for some $A\in\A_0$, then $\B'=\varnothing$.
Similarly, if $B\subseteq K\setminus S$ for some $B\in\B_0$,
then $\A'=\varnothing$. These two events are disjoint,
since $A\cap B\ne\varnothing$ for all $A\in\A_0$ and
$B\in\B_0$.

Suppose first that $a=3$, or that $a=2$ and neither
$\A_0$ nor $\B_0$ consists of two disjoint sets.
Applying Claim~\ref{cl:singleton-small-sets} to the two
ISPs constructed above, we obtain
\[
\begin{aligned}
 \PP\bigl(R(\A',\B')\ne\varnothing\bigr)
 &\le
 1-\PP(A\subseteq S\text{ for some }A\in\A_0)\\
 &\qquad
   -\PP(B\subseteq K\setminus S\text{ for some }B\in\B_0)\\
 &\le
 1-\frac{|\bigcup\A_0|+|\bigcup\B_0|}{2^{2a-1}}\\
 &\le
 1-\frac{|K|+|H|}{2^{2a-1}}.
\end{aligned}
\]
Together with Claim~\ref{cl:singleton-restriction}
and~\eqref{eq:singleton-weight-on-K}, this gives
\[
 \sum_{v\in U}2^{-w(v)}
 \le
 \frac{|K|+|H|}{2^{2a-1}}
 +1-\frac{|K|+|H|}{2^{2a-1}}
 =1.
\]

Finally, suppose that $a=2$ and $\A_0=\{A_1,A_2\}$, where $A_1 \cap A_2 = \emptyset$. In this case for each $B\in\B_0$, $B \subset A_1 \cup A_2$. Hence $|K|=|A_1 \cup A_2|=4$. Minimality of $\A_0$ implies that each of $A_1, A_2$
contains an element of $H$, so $2\le|H|\le4$.

Suppose that $R(\A',\B')\ne\varnothing$. Then there exist
$x\in U\setminus K$ and sets $A\in\A$, $B\in\B$ such that
\[
 A\cap B=\{x\},\qquad
 A\cap K\subseteq S,\qquad B\cap K\subseteq K\setminus S.
\]
Since $B$ intersects the disjoint sets $A_1,A_2$ and contains
$x$, we have $|B|\ge3$.

Fix $v\in H$ and write the set in $\A_0$ containing $v$
as $\{v,u\}$. If $u\notin B$, then $B\cap\{v,u\}=\{v\}$.
This would give
\[
 w(v)\ge|\{v,u\}|+|B|-2=|B|\ge3,
\]
contrary to $w(v)=2$. Hence $u\in B$, and therefore $u\notin S$.

For distinct $v\in H$, the corresponding elements $u$
are distinct, since $A_1$ and $A_2$ are disjoint. Thus $|H|$ fixed elements of $K$ cannot belong to $S$. Note that $A \cap K$ is nonempty. Therefore $S$ must be a nonempty subset of the remaining $4-|H|$ elements.
It follows that
\[
 \PP\bigl(R(\A',\B')\ne\varnothing\bigr)
 \le\frac{2^{4-|H|}-1}{16}.
\]
By Claim~\ref{cl:singleton-restriction}
and~\eqref{eq:singleton-weight-on-K},
\[
 \sum_{v\in U}2^{-w(v)}
 \le\frac{|H|}{4}+\frac{4-|H|}{8}
     +\frac{2^{4-|H|}-1}{16}
 \le1,
\]
where the last inequality holds for $2\le|H|\le4$.
This completes the proof of Theorem~\ref{thm:singleton}.

\section*{Statement on the use of AI}

The authors used ChatGPT (GPT-5.6 and GPT-6 Astra) to assist
in finding and simplifying proofs and in preparing the manuscript. All arguments have been checked by the authors, who take full responsibility for the content.

\end{document}